\documentclass[12pt]{article}
\usepackage{amsmath,amsfonts,latexsym,amsthm,amssymb}
\newcommand{\Qp}{\mathbb Q_p}
\newcommand{\D}{D^\alpha}
\newcommand{\II}{I^\alpha}
\newcommand{\IO}{I^1}
\newcommand{\IOI}{I_0^1}
\newcommand{\DO}{D_O^\alpha}
\newcommand{\HH}{\mathcal H}
\newcommand{\DOI}{D_O^1}

\numberwithin{equation}{section}

\DeclareMathOperator{\const}{const}

\begin{document}

\newtheorem{lem}{Lemma}
\newtheorem{teo}{Theorem}
\newtheorem{cor}{Corollary}
\newtheorem{prop}{Proposition}

\pagestyle{plain}
\title{Non-Archimedean Radial Calculus: Volterra Operator and Laplace Transform}
\author{Anatoly N. Kochubei
\\ \footnotesize Institute of Mathematics,\\
\footnotesize National Academy of Sciences of Ukraine,\\
\footnotesize Tereshchenkivska 3, Kiev, 01024 Ukraine\\
\footnotesize E-mail: \ kochubei@imath.kiev.ua}
\date{}
\maketitle

\vspace*{3cm}
\begin{abstract}
In an earlier paper (A. N. Kochubei, {\it Pacif. J. Math.} 269 (2014), 355--369), the author considered a restriction of Vladimirov's fractional differentiation operator $D^\alpha$, $\alpha >0$, to radial functions on a non-Archimedean field. In particular, it was found to possess such a right inverse $I^\alpha$ that the appropriate change of variables reduces equations with $D^\alpha$ (for radial functions) to integral equations whose properties resemble those of classical Volterra equations. In other words, we found, in the framework of non-Archimedean pseudo-differential operators, a counterpart of ordinary differential equations. In the present paper, we begin an operator-theoretic investigation of the operator $I^\alpha$, and study a related analog of the Laplace transform.
\end{abstract}

\vspace{2cm}


\newpage
\section{Introduction}

The basic linear operator defined on real- or complex-valued functions on a non-Archimedean local field $K$ (such as $K=\Qp$, the field of $p$-adic numbers) is the Vladimirov pseudo-differential operator $\D$, $\alpha >0$, of fractional differentiation \cite{VVZ}; for further development of this subject see \cite{AKS,KKZ,K2001,K2008,Z2016,BK}. Note also the recent publications devoted to applications in geophysical models and to the study of related nonlinear equations \cite{OKM,KK2018,PK,KK2020}.

It was found in \cite{K2014} that properties of $\D$ become much simpler on radial functions. Moreover, in this case it was found to possess a right inverse $\II$, which can be seen as a $p$-adic counterpart of the Riemann-Liouville fractional integral or, for $\alpha =1$, the classical anti-derivative. The change of an unknown function $u=\II v$ reduces the Cauchy problem for an equation with the radial restriction of $\D$ to an integral equation with properties resembling those of classical Volterra equations. In other words, we found, in the framework of non-Archimedean pseudo-differential operators, a counterpart of ordinary differential equations. In \cite{K2014}, we studied linear equations of this kind; nonlinear ones were investigated in \cite{K2020}. Note that radial functions appear as exact solutions of the $p$-adic analog of the classical porous medium equation \cite{KK2018}.

In this paper we study the operator $\IO$ on the ring of integers $O\subset K$ as an object of operator theory. The operator $\IO$ on $L^2(O)$ happens to be a sum of a bounded selfadjoint operator and a simple Volterra operator $\IOI$ with a rank two imaginary part $J$, such that $\operatorname{tr} J=0$. The characteristic matrix-function $W(z)$ of $\IOI$ is such that $W(z^{-1})$ is, in contrast to classical examples, an entire matrix function of zero order.

While the theory of Volterra operators and their characteristic functions is well-developed (see \cite{GK,Br,GT,Ma,Zo}), properties of the above operator are very different from those known for operators of classical analysis and their generalizations. Therefore, while $\IO$ and $\IOI$ are just specific examples, they create a framework for future studies in this area.

Another subject touched in this paper is a version of the Laplace transform. The classical Laplace transform is based on the function $x\mapsto e^{-\lambda x}$ satisfying an obvious differential equation. A similar equation involving $\D$ has a unique radial solution \cite{K2008}. This leads to a definition of the Laplace type transform in the above framework. We prove a uniqueness theorem and the inversion formula for this transform.

\section{Preliminaries}

{\bf 2.1. Local fields.} Let $K$ be a non-Archimedean local field,
that is a non-discrete totally disconnected locally compact
topological field. It is well known that $K$ is isomorphic either
to a finite extension of the field $\mathbb Q_p$ of $p$-adic
numbers (if $K$ has characteristic 0), or to the field of formal
Laurent series with coefficients from a finite field, if $K$ has
a positive characteristic. For a summary of main notions and results
regarding local fields see, for example, \cite{K2001}.

Any local field $K$ is endowed with an absolute value $|\cdot |_K$,
such that $|x|_K=0$ if and only if $x=0$, $|xy|_K=|x|_K\cdot |y|_K$,
$|x+y|_K\le \max (|x|_K,|y|_K)$. Denote $O=\{ x\in K:\ |x|_K\le 1\}$,
$P=\{ x\in K:\ |x|_K<1\}$. $O$ is a subring of $K$, and $P$ is an ideal
in $O$ containing such an element $\beta$ that
$P=\beta O$. The quotient ring $O/P$ is
actually a finite field; denote by $q$ its cardinality. We will
always assume that the absolute value is
normalized, that is $|\beta |_K=q^{-1}$. The normalized absolute
value takes the values $q^N$, $N\in \mathbb Z$. Note that for $K=\mathbb Q_p$
we have $\beta =p$ and $q=p$; the $p$-adic absolute value is normalized.

The additive group of any local field is self-dual, that is if
$\chi$ is a fixed non-constant complex-valued additive character of
$K$, then any other additive character can be written as
$\chi_a(x)=\chi (ax)$, $x\in K$, for some $a\in K$. Below we assume that $\chi$ is a rank
zero character, that is $\chi (x)\equiv 1$ for $x\in O$, while
there exists such an element $x_0\in K$ that $|x_0|_K=q$ and $\chi
(x_0)\ne 1$.

The above duality is used in the definition of the Fourier
transform over $K$. Denoting by $dx$ the Haar measure on the
additive group of $K$ (normalized in such a way that the measure
of $O$ equals 1) we write
$$
\widetilde{f}(\xi )=\int\limits_K\chi (x\xi )f(x)\,dx,\quad \xi
\in K,
$$
where $f$ is a complex-valued function from $L_1(K)$. As usual, the Fourier
transform $\mathcal F$ can be extended from $L_1(K)\cap L_2(K)$ to a
unitary operator on $L_2(K)$. If $\mathcal F f=\widetilde{f}\in L_1(K)$, we
have the inversion formula
$$
f(x)=\int\limits_K\chi (-x\xi )\widetilde{f}(\xi )\,d\xi .
$$

Working with functions on $K$ and operators upon them we often use standard integration formulas; see \cite{K2001,VVZ}. The simplest of them are as follows:

\begin{equation*}
\int\limits_{|x|_K\le q^n}dx=q^n;\quad \int\limits_{|x|_K=q^n}dx=\left( 1-\frac1q \right)q^n.
\end{equation*}

\begin{equation*}
\int\limits_{|x|_K\le q^n}|x|_K^{\alpha -1}\,dx=\frac{1-q^{-1}}{1-q^{-\alpha }}q^{\alpha n};\quad \text{here and above $n\in \mathbb Z,\alpha >0$}.
\end{equation*}

A function $f:\ K\to \mathbb C$ is said to be locally constant, if there exists such an integer $l$ that for any $x\in K$
$$
f(x+x')=f(x), \quad \text{whenever $|x'|\le q^{-l}$}.
$$

The vector space $\mathcal D(K)$ of all locally constant functions with compact supports is used as a space of test functions in analysis on $K$. Note that the Fourier transform preserves $\mathcal D(K)$. There exists a well-developed theory of distributions on local fields; see \cite{AKS,K2001,VVZ}.

\medskip
{\bf 2.2. Vladimirov's operator}. On a test function $\varphi \in \mathcal D(K)$, the fractional differentiation operator $D^\alpha$, $\alpha >0$, is defined as
\begin{equation}
\label{2.1}
\left( D^\alpha \varphi \right) (x)=\mathcal F^{-1}\left[ |\xi |_K^\alpha
(\mathcal F (\varphi ))(\xi )\right] (x).
\end{equation}
Note that $\D$ does not preserve $\mathcal D(K)$; see \cite{AKS} regarding the spaces of test functions and distributions preserved by this operator.

The operator $D^\alpha$ can also be represented as a hypersingular integral operator:
\begin{equation}
\label{2.2}
\left( D^\alpha \varphi \right) (x)=\frac{1-q^\alpha }{1-q^{-\alpha
-1}}\int\limits_K |y|_K^{-\alpha -1}[\varphi (x-y)-\varphi (x)]\,dy
\end{equation}
\cite{K2001,VVZ}.
In contrast to (\ref{2.1}), the expression in the right of (\ref{2.2}) makes sense for wider classes of functions. In particular, $\D$ is defined on constant functions and annihilates them. Denote for brevity $\theta_\alpha =\dfrac{1-q^\alpha }{1-q^{-\alpha -1}}$.

Below we consider the operator $\D$ on a radial function $u=u(|x|_K)$; here we identify the function $x\mapsto u(|x|_K)$ on $K$ with the function $|x|_K\mapsto u(|x|_K)$ on $q^{\mathbb Z}$. This abuse of notation does not lead to confusion.

The explicit expression of $\D u$ for a radial function $u$ satisfying some growth restrictions near the origin and infinity was found in \cite{K2014}. If $u=u(|x|_K)$ is such that
\begin{equation}
\label{2.3}
\sum\limits_{k=-\infty}^m q^k\left| u(q^k)\right| <\infty ,\quad \sum\limits_{l=m}^\infty q^{-\alpha l}\left| u(q^l)\right| <\infty,
\end{equation}
for some $m\in \mathbb Z$, then for each $n\in \mathbb Z$ the expression in the right-hand side of (\ref{2.2}) with $\varphi (x)=u(|x|_K)$ exists for $|x|_K=q^n$, depends only on $|x|_K$, and
\begin{multline}
\label{2.4}
(D^\alpha u)(q^n)=\theta_\alpha \left(1-\frac1q \right)q^{-(\alpha +1)n}\sum\limits_{k=-\infty}^{n-1} q^ku(q^k) +q^{-\alpha n-1}\frac{q^\alpha +q-2}{1-q^{-\alpha -1}}u(q^n)\\
+\theta_\alpha \left(1-\frac1q \right)\sum\limits_{l=n+1}^\infty q^{-\alpha l}u(q^l).
\end{multline}

Under the conditions (\ref{2.3}), the expression (\ref{2.4}) agrees also with the definition of $\D$ in terms of Bruhat-Schwartz distributions (see Chapter 2 of \cite{VVZ}).

\medskip
{\bf 2.3. The regularized integral}. The fractional integral mentioned in Introduction, was defined in \cite{K2014} initially for $\varphi \in \mathcal D(K)$ as follows:
\begin{equation*}
(I^\alpha \varphi )(x)=(D^{-\alpha}\varphi )(x)-(D^{-\alpha}\varphi )(0) \tag{$*$}
\end{equation*}
where $D^{-\alpha}$ is the right inverse of $\D$ introduced by Vladimirov \cite{VVZ}:
\begin{equation*}
\left( D^{-\alpha} \varphi \right) (x)=(f_\alpha *\varphi )(x)=\frac{1-q^{-\alpha} }{1-q^{\alpha
-1}}\int\limits_K |x-y|_K^{\alpha -1}\varphi (y)\,dy,\quad \alpha \ne 1,
\end{equation*}
\begin{equation*}
\left( D^{-1}\varphi \right) (x)=\frac{1-q}{q\log q}\int\limits_K \log |x-y|_K\varphi (y)\,dy.
\end{equation*}

$D^{-1}$ is a right inverse to $D^1$ only on such functions $\varphi$ that
$$
\int\limits_K\varphi (x)\,dx=0.
$$
On such a function $\varphi$ we have also $D^{-1}D^1\varphi =\varphi$.

The above definition $(*)$ leads to explicit expressions
\begin{equation*}
(I^\alpha \varphi )(x)=\frac{1-q^{-\alpha}}{1-q^{\alpha -1}}\int\limits_{|y|_K\le |x|_K}\left( |x-y|_K^{\alpha -1}-|y|_K^{\alpha -1}\right) \varphi (y)\,dy,\quad \alpha \ne 1,
\end{equation*}
and
\begin{equation*}
(I^1\varphi )(x)=\frac{1-q}{q\log q}\int\limits_{|y|_K\le |x|_K}\left( \log |x-y|_K-\log |y|_K\right) \varphi (y)\,dy.
\end{equation*}
Note that the integrals are taken, for each fixed $x\in K$, over bounded sets, and $(I^\alpha \varphi )(0)=0$. These properties are different from those of the anti-derivatives $D^{-\alpha }$ studied in \cite{VVZ}.

Let $u=u(|x|_K)$ be a radial function, such that
$$
\sum\limits_{k=-\infty}^m \max \left( q^k,q^{\alpha k}\right) \left| u(q^k)\right| <\infty ,\quad \text{if $\alpha \ne 1$},
$$
and
$$
\sum\limits_{k=-\infty}^m |k|q^k \left| u(q^k)\right| <\infty ,\quad \text{if $\alpha =1$},
$$
for some $m\in \mathbb Z$. Then \cite{K2014} $I^\alpha u$ exists, it is a radial function, and for any $x\ne 0$,
\begin{equation*}
(I^\alpha u)(|x|_K)=q^{-\alpha}|x|_K^\alpha u(|x|_K)+\frac{1-q^{-\alpha}}{1-q^{\alpha -1}}\int\limits_{|y|_K< |x|_K}\left( |x|_K^{\alpha -1}-|y|_K^{\alpha -1}\right) u(|y|_K)\,dy,\quad \alpha \ne 1,
\end{equation*}
and
\begin{equation}
\label{2.5}
(I^1 u)(|x|_K)=q^{-1}|x|_K u(|x|_K)+\frac{1-q}{q\log q}\int\limits_{|y|_K<|x|_K}\left( \log |x|_K-\log |y|_K\right) u(|y|_K)\,dy.
\end{equation}
On an appropriate class of radial functions, $\II$ is a right inverse to $\D$ \cite{K2014}.An important difference between $D^{-\alpha }$ and $\II$ is the bounded integration domain in the integral formulas for $\II$.

\bigskip
{\bf 2.4. Radial eigenfunctions of $\D$.}

The operator $\D$ defined initially on $\mathcal D(K)$ is, after its closure in $L^2(K)$, a selfadjoint operator with a pure point spectrum $\{ q^{\alpha N},N\in \mathbb Z\}$ of infinite multiplicity and a single limit point zero.

It was shown in \cite{K2008} that for each $N\in \mathbb Z$, there exists a unique (up to the multiplication by a constant) radial eigenfunction
\begin{equation}
\label{2.6}
v_N(|x|_K)=\begin{cases}
1, & \text{if $|x|_K\le q^{-N}$,}\\
-\frac1{q-1}, & \text{if $|x|_K=q^{-N+1}$,}\\
0, & \text{if $|x|_K\ge q^{-N+2}$,}\end{cases}
\end{equation}
corresponding to the eigenvalue $\lambda=q^{\alpha N}$. Below we interpret this function as an analog of the classical exponential function $x\mapsto e^{-\lambda x}$. Note that $v_N\in \mathcal D(K)$; this is a purely non-Archimedean phenomenon reflecting the unusual topological property of $K$, its total disconnectedness.

The operator $\DO$ in the space $L^2(O)$ on the ring of integers (unit ball) $O$ is defined as follows. Extend a function $\varphi \in \mathcal D(O)$ (that is a function $\varphi \in \mathcal D(K)$ supported in $O$) onto $K$ by zero. Apply $\D$ and consider the resulting function on $O$. After the closure in $L^2(O)$ we obtain a selfadjoint operator $\DO$ with a discrete spectrum \cite{K2001,VVZ} (here we do not touch different definitions from \cite{BGPW} and \cite{K2018}).

Denote by $\HH$ the subspace in $L^2(O)$ consisting of radial functions. The functions $v_N$, $N=1,2,\ldots$ belong to $\HH$, as well as the function
$$
v_0(|x|_K)\equiv 1,\quad |x|_K\le 1.
$$
By the definition of $\DO$, the functions $v_N$ are its eigenfunctions corresponding to the eigenvalues $q^{\alpha N}$. As for $v_0$, it is also an eigenfunction, with the eigenvalue $\mu_0=\dfrac{q-1}{q^{\alpha +1}-1}q^\alpha$ \cite{K2001,VVZ}. Therefore $\{ v_N\}_{N\ge 0}$ is an orthonormal system in $L^2(O)$, hence in $\HH$.

We have $\| v_0\| =1$ ($\| \cdot \|$ is the norm in $\HH$),
$$
\|v_N\|^2=\int\limits_{|x|_K\le q^{-N}}dx+(q-1)^{-2} \int\limits_{|x|_K=q^{-N+1}}dx=q^{-N}+(q-1)^{-2}q^{-N+1}(1-\frac1q)=(q-1)^{-1}q^{-N},
$$
$$
\int\limits_{|x|_K\le 1}v_N(|x|_K)\,dx=0,\quad N\ge 1.
$$
Therefore the functions
\begin{equation}
\label{2.7}
e_0( |x|_K)\equiv 1;\quad e_N( |x|_K)=(q-1)^{1/2}q^{N/2}v_N(|x|_K),\quad N\ge 1,
\end{equation}
form an orthonormal system in $\HH$.

\medskip
\begin{lem}
The system $\{ e_N\}_{N\ge 0}$ is an orthonormal basis in $\HH$.
\end{lem}

\medskip
{\it Proof}. Let $u\in \HH$ be orthogonal to all the functions $e_N$. Then
$$
\int\limits_{|x|_K\le 1}u(|x|_K)\,dx=0,
$$
so that
\begin{equation}
\label{2.8}
\sum\limits_{j=-\infty}^0 u(q^j)q^j=0
\end{equation}
and
$$
\int\limits_{|x|_K\le q^{-N}}u(|x|_K)\,dx -(q-1)^{-1}\int\limits_{|x|_K=q^{-N+1}}u(|x|_K)\,dx=0,
$$
so that
\begin{equation}
\label{2.9}
\sum\limits_{j=-\infty}^{-N} u(q^j)q^j-(q-1)^{-1}u(q^{-N+1})=0,\quad N=1,2,\ldots .
\end{equation}

Subtracting from (\ref{2.8}) the equality (\ref{2.9}) with $N=1$, we find that $u(1)=0$. Now the equality (\ref{2.9}) with $N=1$ takes the form
$$
\sum\limits_{j=-\infty}^{-1} u(q^j)q^j=0,
$$
while (\ref{2.9}) with $N=2$ yields
$$
\sum\limits_{j=-\infty}^{-2} u(q^j)q^j-(q-1)^{-1}u(q^{-1})=0.
$$
Subtracting we obtain that $u(q^{-1})=0$.

Repeating the above reasoning we find that $u=0$.\qquad $\blacksquare$

\medskip
Another (obvious) orthonormal basis in $\HH$ is
\begin{equation}
\label{2.10}
f_n( |x|_K)=\begin{cases}
(1-\frac1q)^{-1/2}q^{n/2}, & \text{if $|x|_K=q^{-n}$;}\\
0, & \text{elsewhere,}\end{cases}
\quad n=0,1,2,\ldots .
\end{equation}

The next result is of some independent interest.

\medskip
\begin{prop}
The set of ``polynomials''
\begin{equation}
\label{2.11}
u(|x|_K)=\sum\limits_{n=1}^N a_n|x|_K^n,\quad a_n\in \mathbb C,\ N\ge 1,
\end{equation}
is dense in $\HH$.
\end{prop}

\medskip
{\it Proof}. Suppose that a function $F\in \HH$ is orthogonal to all the functions $X_l(|x|_K)=|x|_K^l$, $l\ge 1$. Using the basis (\ref{2.10}), write
$$
F=\sum\limits_{n=0}^\infty c_nf_n,\quad \{c_n\} \in l^2.
$$
We have
$$
\langle X_l,f_n\rangle =(1-\frac1q)^{-1/2}q^{n/2}\int\limits_{|x|_K=q^{-n}}|x|_K^l\,dx=(1-\frac1q)^{1/2}q^{-n/2-nl},
$$
so that
$$
\langle F,X_l\rangle = (1-\frac1q)^{1/2} \sum\limits_{n=0}^\infty c_nq^{-n/2-nl}=0,\quad l=1,2,\ldots .
$$

Denoting $\beta =q^{-1}$, $b_n=c_nq^{-n/2}$, we see that the vector $(b_0,b_1,b_2,\ldots )\in l^2$ is orthogonal in $l^2$ to each vector $(1,\beta^l,\beta^{2l},\ldots )$, $l\ge 1$. It is known  (\cite{Ha}, Problem 6) that the set of all these vectors is total in $l^2$, so that $F=0$. \qquad $\blacksquare$

\medskip
In fact, the above reasoning proves the density of polynomials (\ref{2.11}) in a wider weighted space determined by the condition $\{ c_nq^{-n/2}\} \in l^2$.

\section{Integration Operators}

{\bf 3.1. The operator $\IO$.} Let us study $\IO$ as an operator in $\HH$, find its matrix representation with respect to the basis $\{e_N\}$ and investigate the spectrum of $\IO$.

\medskip
\begin{prop}
The operator $\IO$ has the matrix representation
$$
\IO =\begin{pmatrix}
0 & -(q-1)^{1/2}q^{-1/2} & -(q-1)^{1/2}q^{-1} & \ldots & -(q-1)^{1/2}q^{-n/2} & \ldots \\
0 & q^{-1} & 0 & \ldots & 0 & \ldots \\
0 & 0 & q^{-2} & \ldots & 0 & \ldots \\
\hdotsfor{6}\\
0 & 0 & 0 & \ldots & q^{-n} & \ldots \\
\hdotsfor{6}
\end{pmatrix}
$$
(only the first row and the principal diagonal have nonzero elements). $\IO$ is a Hilbert-Schmidt operator. Apart from being a point of essential spectrum, $\lambda =0$ is a simple eigenvalue. In addition, $I^1$ has simple eigenvalues $\lambda_m=q^{-m}$, $m=1,2,\ldots$.
\end{prop}

\medskip
{\it Proof}. Since the integral of each function $e_N$, $N\ge 1$, equals zero, we have $D^{-1}D^1 e_N=e_N$. On the other hand, $D^1e_N=q^Ne_N$, so that $D^{-1}e_N=q^{-N}e_N$, and by the definition $(*)$ of $\IO$,
\begin{equation}
\label{3.1}
\IO e_N=q^{-N}e_N-(q-1)^{1/2}q^{-N/2}e_0,\quad N=1,2,\ldots .
\end{equation}

Next, $(\IO e_0)(|x|_K)$, $|x|_K\le 1$, depends only on the values of $e_0$ for $|x|_K\le 1$. Let $f(x)\equiv 1$, $x\in K$. Then $\IO f=0$ \cite{K2014}, so that
\begin{equation}
\label{3.2}
\IO e_0=0\quad \text{in $\HH$}.
\end{equation}

The equalities (\ref{3.1}) and (\ref{3.2}) imply the required matrix representation, which implies the Hilbert-Schmidt property.

Let us find the eigenvalues of $I^1$. As we have seen, $\IO e_0=0$. Suppose that
$$
u=\sum\limits_{n=0}^\infty c_ne_n,\quad \{c_n\}\in l^2,\quad \IO u=\lambda u.
$$
By (\ref{3.1}) and (\ref{3.2}),
$$
\IO u=\sum\limits_{n=1}^\infty q^{-n}c_ne_n-\left[ \sum\limits_{n=1}^\infty (q-1)^{1/2}q^{-n/2}c_n\right] e_0,
$$
and we find that
\begin{equation}
\label{3.3}
\left\{
\begin{aligned}
\lambda c_0 & =-(q-1)^{1/2}\sum\limits_{n=1}^\infty q^{-n/2}c_n;\\
\lambda c_n & =q^{-n}c_n,\quad n\ge 1.
\end{aligned}
\right.
\end{equation}

A nonzero value of $c_n (n\ge 1)$ is possible only for a single index $n=m$, and in this case $\lambda =q^{-m}$. Then the first equation in (\ref{3.3}) gives $c_m=-(q-1)^{-1/2}q^{-m/2}c_0$, so that
$$
u=c_0e_0-(q-1)^{-1/2}q^{-m/2}c_0e_m
$$
is the unique (up to the multiplication by a constant) eigenfunction. \qquad $\blacksquare$

\bigskip
{\bf 3.2. A local representation.} The definition $(*)$ of the operator $\II$ involves operators in $L^2(K)$; then we make restrictions to $L^2(O)$ and $\HH$. In this section we show, for the case where $\alpha =1$, that a similar representation containing only operators in $L^2(O)$ is also possible.

\medskip
\begin{teo}
If $u\in L^2(O)$, then
\begin{equation}
\label{3.4}
(\IO u)(x)=\left( \left( \DOI \right)^{-1}u\right) (x)-\left( \left( \DOI \right)^{-1}u\right) (0).
\end{equation}
\end{teo}

\medskip
{\it Proof}. In \cite{K2018}, we found the resolvent $\left( \DOI -\mu_0 +\mu\right)^{-1}$ where $\mu_0=\dfrac{q}{q+1}$ (the first eigenvalue of $\DOI$), $\mu >0$. In \cite{K2018}, in connection with nonlinear equations, we considered operators in $L^1(O)$, but the result is valid for $L^2(O)$ too. For $\mu =\mu_0$,
\begin{equation}
\label{3.5}
\left( \DOI \right)^{-1}u (x)=\int\limits_{|\xi |_K\le 1}\mathcal K (x-\xi )u(\xi )\,d\xi+\mu_0^{-1}\int\limits_{|\xi |_K\le 1}u(\xi )\,d\xi ,
\end{equation}
where for $|x|_K=q^m$, $m\le 0$,
$$
\mathcal K(x)=\int\limits_{q\le |\eta |_K\le q^{-m+1}}|\eta |_K^{-1}\chi (\eta x)\,d\eta .
$$

Using the well-known integration formula (see, for example, Section 1.5 in \cite{K2001}), we get
\begin{multline*}
\mathcal K(x)=\sum\limits_{j=1}^{-m+1}q^{-j}\int\limits_{|\eta |_K=q^j}\chi (\eta x)\,d\eta =(1-\frac1q)\sum\limits_{j=1}^{-m}1-q^{-1}\\
=-(1-\frac1q)m-q^{-1}=\frac{1-q}{q\log q}\log |x|_K-q^{-1}.
\end{multline*}

By (\ref{3.5}),
$$
\left( \DOI \right)^{-1}u (x)=\frac{1-q}{q\log q}\int\limits_{|\xi |_K\le 1}\log |x-\xi |_Ku(\xi )\,d\xi +\int\limits_{|\xi |_K\le 1}u(\xi )\,d\xi .
$$
Comparing with the expression for $\IO$ and noticing that $|x-\xi |_K-|\xi|_K=0$, if $|\xi |_K>|x|_K$, we obtain (\ref{3.4}).\qquad $\blacksquare$

\bigskip
{\bf 3.3. The Volterra operator.} Let us consider the integral part of (\ref{2.5}), the operator
$$
\left( \IOI u\right) (x)=\frac{1-q}{q\log q}\int\limits_{|y|_K<|x|_K}\left( \log |x|_K-\log |y|_K\right) u(|y|_K)\,dy.
$$
Recall \cite{GK} that a compact operator is called a Volterra operator, if its spectrum consists of the unique point $\lambda =0$. An operator $A$ is called simple, if $A$ and $A^*$ have no common nontrivial invariant subspace, on which these operators coincide. It is known \cite{GK} that a Volterra operator $A$ is simple, if and only if the equations $Af=0$ and $A^*f=0$ have no common nontrivial solutions.

The main technical tool in the study of $\IOI$ is the identity \cite{K2014}
\begin{equation}
\label{3.6}
\int\limits_{|y|_K<|x|_K}\left( \log |x|_K-\log |y|_K\right)|y|_K^m\,dy=d_m|x|_K^{m+1},\quad m=0,1,2,\ldots ,
\end{equation}
where $0<d_m\le Aq^{-m}$, $A>0$ does not depend on $m$.

\medskip
\begin{teo}
The operator $\IOI$ in $\HH$ is a simple Volterra operator with a rank 2 imaginary part $J=\frac1{2i}(A-A^*)$, such that $\operatorname{tr}J=0$.
\end{teo}

\medskip
{\it Proof}. 1) Suppose that $\IOI u=\lambda u,$ $u\in \HH$, $\lambda \in \mathbb C$, $\lambda \ne 0$. Then for $|x|_K\le 1$,
\begin{multline*}
|u(|x|_K)|\le \frac{c}{|\lambda|}\|u\|_{L^2(O)}\left[ \int\limits_{|y|_K<|x|_K}\left( \log |x|_K-\log |y|_K\right)^2\,dy\right]^{1/2}\\
\le \frac{c}{|\lambda|}\|u\|_{L^2(O)}\left[ q^{-1}|x|_K(\log |x|_K)^2\right]^{1/2}\le H
\end{multline*}
where $c=\dfrac{q-1}{q\log q}$, $H$ is a positive constant.

This implies the estimate
$$
|u(|x|_K)|\le \frac{cH}{|\lambda|}\int\limits_{|y|_K<|x|_K}\left( \log |x|_K-\log |y|_K\right)\,dy,
$$
and by the identity (\ref{3.6}) with $m=0$,
$$
|u(|x|_K)|\le \frac{cHA}{|\lambda|}|x|_K.
$$

Similarly, the identity (\ref{3.6}) with $m=1$ gives
$$
|u(|x|_K)|\le \frac{c^2HA^2}{|\lambda|^2}q^{-1}|x|_K^2,
$$
and we find by induction that
\begin{equation}
\label{3.7}
|u(|x|_K)|\le \frac{c^{m+1}HA^{m+1}}{|\lambda|^{m+1}}q^{-1}q^{-2}\cdots q^{-m+1}|x|_K^{m+1},
\end{equation}
for an arbitrary natural number $m$.

Note that
$$
q^{-1}q^{-2}\cdots q^{-m+1}=\left( \frac1q\right)^{m(m-1)/2}.
$$
Together with (\ref{3.7}), this shows that $u\equiv 0$.

\medskip
2) It follows from the definition of $\IOI$ that $\lambda =0$ is an eigenvalue corresponding to the eigenfunction
\begin{equation}
\label{3.8}
u_0(|x|_K)=\begin{cases}
1, & \text{if $|x|_K=1$};\\
0, & \text{if $|x|_K<1$}.\end{cases}
\end{equation}
Let us show that $\lambda =0$ does not correspond to other eigenfunctions.

Suppose that $\IOI \varphi =0$ for some $\varphi \in \HH$, so that
\begin{equation}
\label{3.9}
\sum\limits_{j=-\infty}^{n-1}(n-j)q^j\varphi (q^j)=0,\quad n=0,-1,-2,\ldots .
\end{equation}
Together with (\ref{3.9}), consider a similar equality with $n-1$ substituted for $n$, that is
\begin{equation}
\label{3.10}
\sum\limits_{j=-\infty}^{n-2}(n-1-j)q^j\varphi (q^j)=0,\quad n=0,-1,-2,\ldots .
\end{equation}
Subtracting (\ref{3.9}) from (\ref{3.10}) we find that
$$
\sum\limits_{j=-\infty}^{n-1}q^j\varphi (q^j)=0,\quad n=0,-1,-2,\ldots ,
$$
that is, in particular,
$$
q^{n-1}\varphi (q^{n-1})+q^{n-2}\varphi (q^{n-2})+\cdots =0,
$$
$$
q^{n-2}\varphi (q^{n-2})+q^{n-3}\varphi (q^{n-3})+\cdots =0,
$$
Subtracting the second equality from the first one, we find that $\varphi (q^{-1})=\varphi (q^{-2})=\ldots =0$, so that $\varphi$ is proportional to the eigenfunction (\ref{3.8}).

\medskip
3) The imaginary part $J$ has the following matrix representation with respect to the basis $\{ e_N\}$:
\begin{equation}
\label{3.11}
J=\frac{(q-1)^{1/2}}{2i}\begin{pmatrix}
0 & -q^{-1/2} & -q^{-1} & \ldots & -q^{-N/2} & \ldots \\
q^{-1/2} & 0 & 0 & \ldots & 0 & \ldots \\
q^{-1} & 0 & 0 & \ldots & 0 & \ldots \\
\hdotsfor{6}\\
q^{-N/2} & 0 & 0 & \ldots & 0 & \ldots \\
\hdotsfor{6}
\end{pmatrix}.
\end{equation}

It is easy to write an integral representation
\begin{equation}
\label{3.12}
(Ju)(|x|_K)=\frac{1-q}{2iq\log q}\int\limits_{|y|_k\le 1}(\log |x|_K-\log |y|_K)u(|y|_K\,dy ,
\end{equation}
that is
\begin{equation}
\label{3.13}
(Ju)(|x|_K)=\frac{1-q}{2iq\log q}\left[ \langle u,1\rangle \log |\cdot |_K-\langle u,\log |\cdot |_K\rangle 1\right]
\end{equation}
hence $J$ is a rank 2 operator. We see from (\ref{3.11}) that $\operatorname{tr}J=0$.

\medskip
4) The only solution in $\HH$ (up to the multiplication by a constant) of the equation $\IOI u=0$ is the eigenfunction $u_0$ given by (\ref{3.8}). Suppose that $\left( \IOI \right)^*u_0=0$. Then $Ju_0=0$. However by (\ref{3.12}),
$$
(Ju_0)(|x|_K)=\frac{1-q}{2iq\log q}\log |x|_K\int\limits_{|y|_K=1}dy=-\frac{(q-1)^2}{2iq^2\log q}\log |x|_K,
$$
so that $Ju_0\not\equiv 0$, and we have come to a contradiction. This proves that $\IOI$ is a simple Volterra operator. \qquad $\blacksquare$

Let us calculate the action of $\IOI$ upon the basis $\{f_n\}$ defined in (\ref{2.10}). We find for $|x|_K=q^{-j},j\ge 0$, that
\begin{multline*}
\left( \IOI f_n\right) (|x|_K)=-\frac{(1-q^{-1})^{1/2}}{\log q}q^{n/2}\int\limits_{|y|_K<q^{-j},|y|_K=q^{-n}}(\log |x|_K-\log |y|_K)\,dy\\
=(1-q^{-1})^{1/2}q^{n/2}(j-n)\int\limits_{|y|_K<q^{-j},|y|_K=q^{-n}}dy
=\begin{cases}
(1-q^{-1})^{3/2}q^{-n/2}(j-n), & \text{if $n>j$;}\\
0, & \text{if $n\le j$}.
\end{cases}
\end{multline*}

This implies the equality
$$
\langle \IOI f_n,f_j\rangle =0 \text{ for $n\le j$},
$$
meaning that $\{ f_n\}$ {\it is a basis of triangular representation for the operator} $\IOI$.

\medskip
{\bf Remark.} The operator $\IOI$ is $S$-real with respect to the involution $S$ in $\HH$ given by the complex conjugation. Therefore it is $S$-unicellular (\cite{GK}, Appendix, Theorem 5.5). It is not clear whether it is unicellular in the usual (complex) sense. However it is unicellular in a smaller space $\HH^p$ defined as a completion of the set of all ``polynomials'' $\varphi (|x|_K)=\sum\limits_{j=0}^N c_j|x|_K^j$ with respect to the norm $\|\varphi \|=\left\{ \sum |c_j|^p\right\}^{1/p}$,$1\le p<\infty$. By virtue of (\ref{3.6}), $\IOI$ acts on the space $\HH^p$ (isomorphic to $l^p$) as a weighted shift, for which the unicellularity was proved by Yakubovich \cite{Ya}.

\bigskip
{\bf 3.4. Characteristic function.} Following the notation in \cite{GT}, let us write (\ref{3.13}) in the form
$$
\frac1i \left( \IOI-\left( \IOI\right)^*\right) u=\sum\limits_{\alpha,\beta =1}^2\langle u,h_\alpha \rangle j_{\alpha \beta}h_\beta
$$
where $h_1(|x|_K)=\dfrac{q-1}{iq\log q}(=\const)$, $h_2(|x|_K)=-\log |x|_K$, $x\in O$, $j=\left( \begin{smallmatrix}
0 & 1\\
1 & 0 \end{smallmatrix} \right)$. For the operator $\IOI$, we consider the $2\times 2$ characteristic matrix-function of inverse argument
$$
W(z^{-1})=E+izj\left[ \langle \left( E-z\IOI\right)^{-1}h_\alpha ,h_\beta\rangle \right]_{\alpha,\beta =1}^2
$$
where $E$ denotes both the unit operator in $\HH$ and the unit matrix.

For the Volterra operator $\IOI$, $W(z^{-1})$ is an entire matrix-function.

\medskip
\begin{teo}
Matrix elements of $W(z^{-1})$ are entire functions of zero order.
\end{teo}

\medskip
{\it Proof}. For small values of $|z|$, the Fredholm resolvent $\left( E-z\IOI\right)^{-1}$ is given by the Neumann series
$$
\left( E-z\IOI\right)^{-1}f=\sum\limits_{n=0}^\infty \left( z\IOI \right)^nf,\quad f\in \HH.
$$
In order to calculate the characteristic function,we have to compute the functions $\left(\IOI\right)^n 1$ and $\left(\IOI\right)^n\log |\cdot |_K$. The first of them is obtained easily from (\ref{3.6}):
$$
\left( \left( \IOI \right)^n 1\right)(|x|_K)=c^n\prod\limits_{m=0}^{n-1}d_m\cdot |x|_K^n,\quad |x|_K\le 1,
$$
where $c=\dfrac{1-q}{q\log q}$, $0<d_m\le Aq^{-m}$. Summing the progression we find that
\begin{equation}
\label{3.14}
\left( \left( E-z\IOI\right)^{-1}1\right) (|x|_K)=\sum\limits_{n=0}^\infty\rho_nz^n|x|_K^n,\quad |\rho_n|\le C^nq^{-n^2/2},
\end{equation}
where $C>0$ is a constant.

Let us consider $\left(\IOI\right)^n\log |\cdot |_K$. We have
$$
\left( \IOI \log|\cdot |_K\right) (|x|_K)=c\int\limits_{|y|_K<|x|_K}\left( \log |x|_K-\log |y|_K\right) \log |y|_K\,dy.
$$
Setting $y=xt$, $|t|_K<1$, we obtain
\begin{multline*}
\left( \IOI \log|\cdot |_K\right) (|x|_K)=-c|x|_K\int\limits_{|t|_K<1}\log |t|_K(\log |x|_K+\log |t|_K)\,dt=-ca_0|x|_K\log |x|_k-cb_0|x|_K\\
\stackrel{\text{def}}{=}\sigma_1|x|_K\log |x|_K-\eta_1|x|_K
\end{multline*}
where
$$
a_0=\int\limits_{|t|_K<1}\log |t|_K\,dt,\quad b_0=\int\limits_{|t|_K<1}\log^2 |t|_K\,dt.
$$

A similar calculation yields the expression
$$
\left( \IOI (|\cdot |_K\log |\cdot |_K)\right) (|x|_K)=-ca_1|x|^2_K\log |x|_k-cb_1|x|_K^2
$$
where
$$
a_1=\int\limits_{|t|_K<1}|t|_K\log |t|_K\,dt,\quad b_1=\int\limits_{|t|_K<1}|t|_K\log^2 |t|_K\,dt.
$$
Together with (\ref{3.6}), this implies the formula
$$
\left( \left(\IOI\right)^2 \log|\cdot |_K\right) (|x|_K)=c^2a_0a_1|x|^2_K\log |x|_k+c^2a_0b_1|x|_K^2-cb_0d_1|x|_K^2\stackrel{\text{def}}{=}\sigma_2|x|_K^2\log |x|_K+\eta_2|x|_K^2.
$$

Introducing similar constants for the next iterations,
$$
a_n=\int\limits_{|t|_K<1}|t|^n_K\log |t|_K\,dt,\quad b_n=\int\limits_{|t|_K<1}|t|^n_K\log^2 |t|_K\,dt,
$$
and noticing that $|a_n|,|b_n|\le Mq^{-n}$, we prove by induction that
\begin{equation}
\label{3.15}
\left(\IOI\right)^n\log |\cdot |_K=\sigma_n|x|_K^n\log |x|_K+\eta_n|x|_K^n
\end{equation}
where $|\sigma_n|,|\eta |_n\le C^nq^{-1}q^{-2}\cdots q^{-n+1}=C^nq^{-n(n-1)/2}$.

It follows from (\ref{3.15}) that
$$
\left( \left( E-z\IOI\right)^{-1}\log |\cdot |_K\right) (|x|_K)=\sum\limits_{n=0}^\infty\sigma_nz^n|x|_K^n\log |x|_K+\sum\limits_{n=0}^\infty\eta_nz^n|x|_K^n
$$
where $|\sigma_n|,|\eta_n|\le C_1^nq^{-n^2/2}$.

Now we can compute the matrix-function $W(z^{-1})$. By (\ref{3.14}),
$$
\langle \left( E-z\IOI\right)^{-1}h_1,h_1\rangle =\const \cdot \sum\limits_{n=0}^\infty \rho_nz^n\int\limits_{|x|_K^n\le 1}\,dx=\sum\limits_{n=0}^\infty \gamma_nz^n
$$
where $|\gamma_n|\le C_2^nq^{-n^2/2}$, so that this matrix element is an entire function of zero order. Other matrix elements are estimated similarly on the basis of (\ref{3.15}), by inserting 1 as an upper bound of $|x|_K$ and taking into account the convergence of the integrals of $\log |x|_K$ and $\log^2 |x|_K$. \qquad $\blacksquare$

\section{The Laplace Type Transform}

{\bf 4.1. Definition and Properties.} Our definition of a Laplace type transform is based on the function $v_N$ given by (\ref{2.6}). It is essential that $v_N\in \mathcal D(K)$. As we know, $\D v_N=q^{\alpha N}v_N$ ($\alpha >0$).

Let $\xi \in K$, $|\xi |_K=q^N$. Then for any $x\in K$, $v_N(|x|_K)=v_0(|x\xi |_K)$,
$$
D_x^\alpha v_0(|x\xi |_K)=D_x^\alpha v_N(|x|_K)=q^{\alpha N}v_N(|x|_K)=|\xi|_K^\alpha v_0(|x\xi |_K).
$$

We call the function
$$
\widetilde{\varphi}(|\xi|_K)=\int\limits_Kv_0(|x\xi |_K)\varphi (|x|_K)\,dx
$$
the Laplace type transform of a radial function $\varphi\in L^1_{\text{loc}}(K)$. By the dominated convergence theorem, $\widetilde{\varphi}$ is continuous, bounded, and $\widetilde{\varphi}(|\xi|_K)\to 0$, $|\xi|_K\to \infty$.

As a simple computation shows, if $\varphi (|x|_K)\equiv \const$, then $\widetilde{\varphi}(|\xi|_K)\equiv 0$.

The above calculations, together with the selfadjointness of $\D$ in $L^2(K)$, show that
$$
\widetilde{\D \varphi} (|\xi|_K)=\int\limits_K\left(  D_x^\alpha v_0(|x\xi |_K)\right) (|x|_K)\varphi (|x|_K)\,dx=|\xi|_K^\alpha\widetilde{\varphi}(|\xi|_K),\quad \xi\in K.
$$

\medskip
\begin{teo}[uniqueness]
If $\widetilde{\varphi}(|\xi|_K)\equiv 0$, then $\varphi (|x|_K)\equiv \const$.
\end{teo}

\medskip
{\it Proof.} By the definition,
$$
\widetilde{\varphi}(|\xi|_K)=\int\limits_{|x|_K\le |\xi|_K^{-1}}\varphi (|x|_K)\,dx-\frac1{q-1}\int\limits_{|x|_K=q|\xi|_K^{-1}}\varphi (|x|_K)\,dx.
$$

Let $|\xi|_K=q^n$, $n\in \mathbb Z$. Then
$$
\widetilde{\varphi}(q^n)=(1-\frac1q)\sum\limits_{j=-\infty}^{-n}\varphi (q^j)q^j-\varphi (q^{-n+1})q^{-n}.
$$
If we denote $\widetilde{\varphi}(q^n)=f_n$, then
$$
f_{n+1}=(1-\frac1q)\sum\limits_{j=-\infty}^{-n-1}\varphi (q^j)q^j-\varphi (q^{-n})q^{-n-1},
$$
so that
\begin{equation}
\label{4.1}
f_n-f_{n+1}=q^{-n}\left[ \varphi (q^{-n})-\varphi (q^{-n+1})\right].
\end{equation}

If $\widetilde{\varphi}(q^n)=0$ for all $n$, then, by (\ref{4.1}), $ \varphi (q^{-n})=\varphi (q^{-n+1})$ for all $n$, so that \linebreak $\varphi (q^{-n})\equiv \const$. \qquad $\blacksquare$

\medskip
The identity (\ref{4.1}) is of some independent interest, and we formulate it as a corollary.

\medskip
\begin{cor}
For all $n\in \mathbb Z$,
\begin{equation}
\label{4.2}
\widetilde{\varphi}(q^{n})-\widetilde{\varphi}(q^{n+1})=q^{-n}\left[ \varphi (q^{-n})-\varphi (q^{-n+1})\right].
\end{equation}
\end{cor}

\medskip
\begin{cor}
A function $\varphi$ is (strictly) monotone, if and only if $\widetilde{\varphi}$ is (strictly) monotone.
\end{cor}

\bigskip
{\bf 4.2. Inversion formula.}

\begin{teo}
For each $n=1,2,\ldots $,
\begin{equation}
\label{4.3}
\varphi (q^m)=\varphi (1)+\sum\limits_{j=0}^{m-1}q^{-j}\left[ \widetilde{\varphi}(q^{-j+1})-\widetilde{\varphi}(q^{-j})\right],
\end{equation}
\begin{equation}
\label{4.4}
\varphi (q^{-m})=\varphi (1)+\sum\limits_{j=1}^m q^j\left[ \widetilde{\varphi}(q^j)-\widetilde{\varphi}(q^{j+1})\right],
\end{equation}
\end{teo}

\medskip
{\it Proof}. According to (\ref{4.2}),
\begin{gather*}
\varphi (1)-\varphi (q)=\widetilde{\varphi}(1)-\widetilde{\varphi}(q),\\
\varphi (q)-\varphi (q^2)=q^{-1}\left[ \widetilde{\varphi}(q^{-1})-\widetilde{\varphi}(1)\right],\\
\varphi (q^2)-\varphi (q^3)=q^{-2}\left[ \widetilde{\varphi}(q^{-2})-\widetilde{\varphi}(q^{-1})\right],
\end{gather*}
etc. Summing up the first $m$ equalities we obtain (\ref{4.3}).

Similarly, by (\ref{4.2}),
\begin{gather*}
\varphi (q^{-1})-\varphi (1)=q\left[ \widetilde{\varphi}(q)-\widetilde{\varphi}(q^2)\right],\\
\varphi (q^{-2})-\varphi (q^{-1})=q^2\left[ \widetilde{\varphi}(q^2)-\widetilde{\varphi}(q^3)\right],
\end{gather*}
etc, and the summation yields (\ref{4.4}). \qquad $\blacksquare$

\section*{Acknowledgement}
This work was funded in part under the research project "Markov evolutions in real and $p$-adic spaces" of the Dragomanov National Pedagogic  University of Ukraine.

\end{document}